\documentstyle[11pt]{article}
\newtheorem{theorem}{Theorem}[section]

\newtheorem{lemma}{Lemma}[section]

\newcommand{\U}{U_q(C_n^{(1)})}
\renewcommand{\a}{\alpha}

\newcommand{\la}{\lambda}
\newcommand{\be}{\beta}

\newcommand{\ep}{\epsilon}
\newcommand{\vep}{\varepsilon}
\newcommand{\ta}{\tilde{\alpha}}
\newcommand{\tla}{\tilde{\lambda}}

\title{Level One Representations of Quantum Affine Algebras $U_q(C^{(1)}_n)$}

\author{Naihuan Jing$^1$, Yoshitaka Koyama$^2$\\
 and Kailash C. Misra$^1$}
\begin{document}           
\date{January 7, 1998}
\maketitle                 

\begin{center}
$^1$Department of Mathematics\\
North Carolina State University\\
Raleigh, NC 27695-8205, USA
\end{center}
\begin{center}
$^2$Research Institute for Mathematical Sciences\\
Kyoto University,
Kyoto 606, Japan
\end{center}

\begin{abstract}
We give explicit constructions of quantum symplectic affine algebras
at level $1$ using vertex operators.
\end{abstract}

\section{Introduction}

Vertex operator constructions of affine Lie algebras started with
the work of Lepowsky-Wilson \cite{kn:LW} in the principal picture.
The homogeneous construction was given by Frenkel-Kac \cite{kn:FK}
and Segal \cite{kn:S} for the simply laced types. The non-simply laced
cases were provided in \cite{kn:BT} and \cite{kn:GNOS}. These constructions
led to many applications in mathematics and physics.

The quantization of these constructions for the quantum affine algebras
started with the work of Frenkel-Jing \cite{kn:FJ}, where the simply laced
types were constructed. Subsequently
the twisted types $(ADE)^{(r)}$ \cite{kn:J1},
type $B_n^{(1)}$ \cite{kn:Br, kn:JM} and $G_2^{(1)}$ \cite{kn:J2}
 at level $1$ were constructed.
 Very recently we have provided  explicit
realizations of the quantum symplectic affine algebra $U_q(C_n^{(1)})$
at level $-1/2$ in \cite{kn:JKM},
 which are the so-called admissible representations (with rational levels).
The bosonic realizations of $\U$ at level $-1/2$
were obtained by implicitly quantizing
some $\beta \gamma$-system. All these constructions have been used in
obtaining $q$-vertex operators, which in turn provide solutions to associated
quantum
Knizhnik-Zamolodchikov equations.

In this paper we give an explicit construction of $\U$ at level $1$.
Our constructions can be viewed as a quantization of the constructions
in \cite{kn:BT} and \cite{kn:GNOS}. We consider some auxiliary bosonic
fields to build the Fock space. The action of the quantum affine algebra
on the Fock space is given via vertex operators. We hope that our paper also
further clarifies and simplifies the classical constructions.

The idea is to represent the vertex operator as a sum of two fields, where
each one resembles to the original field operator. This idea was already
present in the classical cases, but one needs new techniques to
be able to quantize the construction.

The realization of level $-1/2$ is easier than the level $1$ case due to
the simpler character formulas and $\beta\gamma$-system. As we saw in
\cite{kn:JKM} the level $-1/2$ construction is also completely free from
cocycle consideration, which is critical in almost all
level one cases.

It turns out that there is really a
 subtlety about the sign changes in the construction. Besides
the usual cocycle construction, we have to incorporate further sign factors
$(-1)^{(1-\ep)\partial_i}$ ($\ep=\pm 1$) to simplify the original
construction of \cite{kn:BT}. It is fair to say that the correct
construction depends on whether we completely understand the behavior
of the sign factors, which may partly explain why our simple construction
is found now.

In all rigorous constructions of quantum affine algebras, the Serre relations
are always the most complicated relations to prove. As explained by Jing in
\cite{kn:J1} and more generally in \cite{kn:J5} one has to be able to prove
some combinatorial identities (cf. (\ref{identity1}, \ref{identity2})).
We use the ideas of
\cite{kn:J1} and \cite{kn:J5} to prove all Serre relations.

The paper is organized as follows.
In section two we give the basic notations and preliminaries. In the next
section the main results are presented. The proof of the main theorem
occupies section four. Our representation is a reducible representation, and
we show that our space contains submodules with all level one dominant
highest weights for $U_q(C_n^{(1)})$.

\section{Quantum affine algebras $U_q(C^{(1)}_n)$}

Let $\a_i=e_i-e_{i+1}$ ($i=1, \cdots, n-1$) and $\a_n=2e_n$
be the simple roots of the simple Lie algebra ${sp}_{2n}$, and $
\lambda_i=e_1+\cdots+e_i$ ($i=1, \cdots, n$) be the fundamental
weights. Let $P={\bf Z}e_1+\cdots+{\bf Z}e_n$ and
$Q={\bf Z}\a_1+\cdots+{\bf Z}\a_n$ be the weight and root lattices.
We then let $\Lambda_i, i=0, \cdots, n$ be the fundamental
weights for the affine Lie algebra $\widehat{sp}_{2n}$, here $\Lambda_i=
\lambda_i+\Lambda_0$. We still use $\a_i$ ($i=1, \cdots, n$) together
with $\a_0$ as the
simple roots for the affine Lie algebra $\widehat{sp}_{2n}$.
The nondegenerate symmetric bilinear form $(\ |\ )$ on ${\bf h}^*$, the
dual Cartan subalgebra of $\widehat{sp}_{2n}$,
satisfies that
\begin{equation}
(e_i|e_j)=\frac12\delta_{ij}, \ \ (\delta|\alpha_i)=(\delta|\delta)=0
\ \ \mbox{for all} \ i,j\in I \label{E:2.2},
\end{equation}
and then the generalized Cartan matrix $(A_{ij})$ of $C_n^{(1)}$ is
given by
$
A_{ij}=\frac{(\a_i|\a_j)}{d_i},
$
where $d_i=\frac{(\a_i|\a_i)}2$ and
$(d_0, \cdots, d_n)=(1, 1/2, \cdots, 1/2, 1)$.

For $q$ generic (not a root of unity) let $q_i=q^{d_i}, i\in I$.
The quantum affine algebra $U_q(C_n^{(1)})$ \cite{kn:Dr} is
the associative algebra with 1 over ${\bf C}(q^{1/2})$
generated by the elements $x_{ik}^{\pm}$, $a_{il}$, $K_i^{\pm 1}$,
$\gamma^{\pm 1/2}$, $q^{\pm d}$ $(i=1,2,\cdots,n, k\in {\bf Z},
l\in {\bf Z} \setminus \{0\})$ with the following defining relations :
\begin{eqnarray}
 [\gamma^{\pm 1/2}, u]&=&0 \ \ \mbox{for all} \ u\in U_q(C_n^{(1)}),
\label{E:R1}\\
\mbox{} [a_{ik}, a_{jl}]&=&\delta_{k+l,0}
\displaystyle\frac {[(\a_i|\a_j)k]}{k}
\displaystyle\frac {\gamma^k-\gamma^{-k}}{q-q^{-1}},
\label{E:R2}\\
\mbox{} [a_{ik}, K_j^{\pm 1}]&=&[q^{\pm d}, K_j^{\pm 1}]=0,
\label{E:R3}\\
 q^d x_{ik}^{\pm} q^{-d}&=&q^k x_{ik}^{\pm }, \ \
q^d a_{il} q^{-d}=q^l a_{il},
\label{E:R4}\\
 K_i x_{jk}^{\pm} K_i^{-1}&=&q^{\pm (\alpha_i|\alpha_j)} x_{jk} ^{\pm},
\label{E:R5}\\
\mbox{} [a_{ik}, x_{jl}^{\pm}]&=&\pm \displaystyle\frac {[(\a_i|\a_j)k]}{k}
\gamma^{\mp |k|/2} x_{j,k+l}^{\pm},
\label{E:R6}\\
(z-q^{\pm(\a_i|\a_j)}w)X_i^{\pm}(z)X_j^{\pm}(w)
           &+&(w-q^{\pm(\a_i|\a_j)}z)X_j^{\pm}(w)X_i^{\pm}(z) =0
\label{E:R7}\\
\mbox{}
[X_i^+(z),X_j^-(w)]
             &=&
\frac{\delta_{ij}}{(q_i-q_i^{-1})zw}
                 \left(
                 \psi_i(w\gamma^{1/2})
                 \delta(\frac{w\gamma}{z})\right.  \nonumber\\
                 &&\qquad -\left.\varphi_i(w\gamma^{-1/2})
                  \delta(\frac{w\gamma^{-1}}z)
                 \right)
\label{E:R8}
\end{eqnarray}
where
$X_i^{\pm}(z)=\sum_{n\in {\bf Z}}x^{\pm}_{i,n}z^{-n-1}$, $\psi_{im}$ and
$\varphi_{im}$ $(m\in {\bf Z}_{\ge 0})$
are defined by
\begin{eqnarray}
&&\psi_i(z)=\sum_{m=0}^{\infty} \psi_{im} z^{-m}
=K_i \textstyle {exp} \left( (q-q^{-1}) \sum_{k=1}^{\infty} a_{ik}
z^{-k}
\right),
\nonumber\\
&&\varphi_i(z)=\sum_{m=0}^{\infty} \varphi_{i,-m} z^{m}
=K_i^{-1} \textstyle {exp} \left(- (q-q^{-1}) \sum_{k=1}^{\infty}
a_{i,-k} z^{k}
\right),
\nonumber\\
&&\mbox{Sym}_{z_1,\cdots,z_m}\sum_{r=0}^{m=1-A_{ij}}(-1)^r
\left[\begin{array}{c} m\\r\end{array}\right]_iX_i^{\pm}(z_1)\cdots
X^{\pm}_{i}(z_r) \nonumber\\
&&\qquad\qquad\qquad \times X^{\pm}_{j}(w)X^{\pm}_{i}(z_{r+1})\cdots
X^{\pm}_{i}
(z_m)=0, \mbox{for}\ i\neq j.
\label{E:R9}
\end{eqnarray}
where $[m]_i=(q_i^{m}-q_i^{-m})/(q_i-q_i^{-1})$.

\section{The main results}

Let $m\in {\bf Z}$. Let
$a_i(m)$ ($i=1, \cdots, n$) and $b_i(m)$ ($i=1, \cdots, n-1$) be
two sets of independent operators satisfying the following
Heisenberg relations:
\begin{equation}
\begin{array}{rcl}
\ [a_i(m), a_j(l)]&=&\delta_{m+l, 0}\frac{[(\a_i|\a_j)]}{m}[m],\\
\ [b_i(m), b_j(l)]&=&\delta_{m+l, 0}\frac{[(\a_i|\a_j)]}{m}[m],\\
\ [a_i(m), b_j(l)]&=&0, \qquad [a_i(0), a_j(m)]=[b_i(0), b_j(m)]=0
\end{array}
\end{equation}

In order to construct the Fock space we introduce an identical copy of
the root lattice of $A_{n-1}$
as the sublattice $\tilde{Q}=Q[A_{n-1}]$ of short roots of $Q$.
The basis of $\tilde{Q}$
will be denoted by $\tilde{\a}_i, i=1, \cdots, n-1$. Thus
$$(\ta_i|\ta_j)=(\a_i|\a_j)=\delta_{ij}-\frac12\delta_{|i-j|, 1}.
$$
We also consider the associated weight lattice $\tilde{P}=P[A_{n-1}]$
defined by the inner product.

The Fock space ${\cal V}$ is defined to be the tensor product
of the symmetric algebra generated by $a_i(-m), b_i(-m)$ and the group
algebra generated by $e^{\la}\otimes e^{\tilde{\la}}$
such that $(\ta_i|\la)\pm(\ta_i|\tilde{\la})\in {\bf Z}$ for each
$i\in\{ 1, \cdots, n\}$, where
$\la\in P$ and $\tilde{\la}\in \tilde{P}$. Note that we treat
$\ta_n=0$.

The action of $a_i(m)$ and $b_i(m)$ with $(m\neq 0)$
on $\cal V$ is obtained by considering
the Fock space $\cal V$ as some quotient space of the Heisenberg algebra
tensored with the group algebras of $P$ and $\tilde{P}$. The operators
$a_i(0)$, $b_i(0)$, $e^{\a}$, $e^{\ta}$ act on $\cal V$ by
\begin{eqnarray*}
a_i(0)e^{\la}e^{\tilde{\la}}&=&(\a_i|\la)e^{\la}e^{\tilde{\la}},\quad
b_i(0)e^{\la}e^{\tilde{\la}}=(\ta_i|\tla)e^{\la}e^{\tilde{\la}},\\
e^{\a}e^{\la}e^{\tilde{\la}}&=&e^{\a+\la}e^{\tilde{\la}},\quad
e^{\ta}e^{\la}e^{\tilde{\la}}=e^{\la}e^{\ta+\tla}.
\end{eqnarray*}

 The normal product $: \quad :$
is defined as usual:
   $$ :a_i(m) a_j(l): = a_i(m) a_j(l) \, (m \leq l) ,
   \; \mbox{or} \ a_j(l) a_i(m) \, (m>l) , $$
   $$ :e^{\a} a_i(0):=:a_i(0) e^{\a}:=e^{\a} a_i(0) \, , $$
   $$ :e^{\ta} b_i(0):=:b_i(0) e^{\ta}:=e^{\ta} b_i(0) \, . $$
   and similarly for product involving the $b_i(m)$.

  The degree operator ${d}$ is defined by
 $$ d. v
      =(m_1+\cdots +m_s+l_1+\cdots +l_t+
(\la|\la)+(\tla|\tla))
       v, $$
 where $v=a_{i_1}(m_1)\cdots a_{i_s}(m_s)b_{j_1}(l_1)\cdots b_{j_t}(l_t)
e^{\la}e^{\tla}$ is a basis element in $\cal V$.

It is easy to see that
$a_i(m), b_j(l), e^{\a_i}, e^{\ta_j}$ commute with each other except that
$$
[a_i(0), e^{\a_j}]=(\alpha_i|\alpha_j)e^{\a_j}
, \qquad [b_i(0), e^{\ta_j}]=(\ta_i|\ta_j)e^{\ta_j}.
$$

Let $\vep( \ \ ,\ \ )$: $P\times P\longrightarrow {\pm 1}$
 be the quasi-cocycle such that
\begin{eqnarray*}
\vep(\a, \be+\theta)&=&\vep(\a,\be)\vep(\a, \theta),\\
\vep(\a+\be, \theta)&=&\vep(\a, \theta)\vep(\be, \theta)
(-1)^{(\overline{\a+\be}-\overline{\a}-\overline{\be}|\overline{\theta})},\\
\vep(\a, \be)\vep(\be, \a)&=&(-1)^{(\a|\be)+(\overline{\a}|\overline{\be})},\\
\vep(\a, \be+\theta)\vep(\be, \theta)&=&
\vep(\a, \be)\vep(\a+\be,
\theta)(-1)^{(\overline{\a+\be}-\overline{\a}-\overline{\be}|\overline{\theta})}.
\end{eqnarray*}
where the $-$ is the projection from $P$ to $\tilde P$ defined by
$$\a=\sum_{i=1}^n m_i\a_i\in P \longmapsto \overline{\a_i}=
\sum_{i=1}^{n-1} \overline{m_i}\a_i, \quad \overline{m_i}=m_i(mod\, 2).
$$

We construct such a cocycle directly by
\begin{equation}
\vep(\a_i,\a_j)=\left\{\begin{array}{ll}
                -1, & i=j\\
                1,  & i<j\\
                (-1)^{A_{ij}}, & i>j \end{array}\right.
\end{equation}
and it is easy to verify that the quasi-cocycle satisfies all the defining
relations. In particular, we have
\begin{equation}
\vep(\a_i,\a_j)\vep(\a_j,\a_i)=\left\{\begin{array}{ll}
                (-1)^{2(\a_i|\a_j)}, & \mbox{if}\ \ 1\leq i, j\leq n-1\\
                (-1)^{(\a_i|\a_j)},  & \mbox{otherwise}
                \end{array}\right.
\end{equation}

For $\a\in P$ we define the operators
$\vep_{\a}$ on $\cal V$ such that
\begin{equation}
\vep_{\a}e^{\la}e^{\tla}=\vep(\a, \la)e^{\la}e^{\tla}.
\end{equation}
then
\begin{equation}
\vep_{\a}\vep_{\be}=\vep(\a, \be)\vep(\be, \a)\vep_{\be}\vep_{\a}.
\end{equation}

We denote $\vep_i=\vep_{\a_i}$ for $i=1, \cdots, n$. In particular, we
have:
\begin{eqnarray*}
\vep_i\vep_{j}&=&(-1)^{2(\a_i|\a_j)}\vep_j\vep_i, \qquad 1\leq i, j\leq n-1,\\
\vep_i\vep_{n}&=&(-1)^{(\a_i|\a_n)}\vep_n\vep_i, \qquad i=1, \cdots, n
\end{eqnarray*}

We can now introduce the main vertex operators.
\begin{eqnarray*}
&Y_i^{\pm}(z)&=
         \exp ( \pm \sum^{\infty}_{k=1}
                 \frac{a_i(-k)}{[k]} q^{\mp \frac{k}{2}} z^k)\\
    &&\qquad     \times\exp ( \mp \sum^{\infty}_{k=1}
           \frac{a_i(k)   }{[k]} q^{\mp \frac{k}{2}} z^{-k})
         e^{\pm \a_i} z^{\pm a_i(0)}\vep_i ,\\
&U_j(z)&=
         \exp ( \sum^{\infty}_{k=1}
                 \frac{b_j(-k)}{[k]}z^k)
         \exp ( - \sum^{\infty}_{k=1}
           \frac{b_j(k)   }{[k]}z^{-k})
         e^{ \ta_j} z^{ b_j(0)},\\
&U_j^*(z)&=
         \exp (-\sum^{\infty}_{k=1}
                 \frac{b_j(-k)}{[k]}z^k)
         \exp (\sum^{\infty}_{k=1}
           \frac{b_j(k)   }{[k]}  z^{-k})
         e^{-\ta_j} z^{- b_j(0)},\\
&Z^{\pm}_j(z)&=U_j(q^{\pm 1/2}z)+(-1)^{2a_j(0)}U_j^*(q^{\mp 1/2}z),
\end{eqnarray*}
where $i\in\{1, \cdots, n\}, j\in\{1, \cdots, n-1\}$. For simplicity
we define $Z_n^{\pm}(z)=1$.

\medskip
\begin{theorem} \label{T:1}The Fock space
    $\cal V$ is a $U_q$-module of level $1$
    under the action defined by
\begin{eqnarray*}
       K_i &\longmapsto& q^{a_i(0)}, \qquad
   q^d \longmapsto q^{{d}}, \\
 a_{im} &\longmapsto& a_i(m),\qquad \gamma \longmapsto q, \\
X_i^{\pm}(z) &\longmapsto&
             Y_i^{\pm}(z)Z_i^{\pm}(z),  \qquad i=1, \cdots n.
\end{eqnarray*}
The module $\cal V$ contains submodules generated by the highest weight vectors
$e^{\la_i}e^{\tla_i}$ with weight $\Lambda_i$ where $i=0, \cdots, n$. Here we
denote
$\la_0=\tla_0=0$.
\end{theorem}

\section{Proof of the main theorem}

In this section we prove theorem \ref{T:1} by vertex operator techniques.
Since the operator $a_i(m)$ commutes with $b_j(m)$ (or the operator
$Y^{\pm}_i(z)$ commute with $U_j(w)$ and  $U_j^*(w)$)
 it is clear that
relations
(\ref{E:R1}-\ref{E:R6}) satisfy the Drinfeld relations.

Following \cite{kn:J3} we use basic hypergeometric
series to simplify the presentation.

For $a\in {\bf R}$ we define
$$(1-z)_{q^2}^a:=\frac{(zq^{-a+1};q^2)_{\infty}}{(zq^{a+1};q^2)_{\infty}}
=exp(-\sum_{n\geq 1}\frac{[an]}{n[n]}z^n)
$$
where $(w;q^2)_{\infty}=\prod_{n=0}^{\infty}(1-wq^{2n})$ is the usual
$q$-number. Note that these $q$-series are defined as power series
in $w$. The first few examples of $q$-polynomials are listed in the
following.
\begin{eqnarray*}
(1-z)_{q^2}^{\pm 1}&=&(1-z)^{\pm 1}, \\
(1-z)_{q^2}^{2}&=&(1-qz)(1-q^{-1}z), \\
(1-z)_{q^2}^{1/2}&=&\frac{(zq^{1/2}; q^2)_{\infty}}{(zq^{3/2};q^2)_{\infty}}.
\end{eqnarray*}

It is clear that
$$(1-z)_{q^2}^{-a}=\frac 1{(1-z)_{q^2}^a}.$$

For $\ep, \ep'\in \{\pm=\pm 1\}$,
the operator product expansions (OPE) are given by:
\begin{eqnarray}
Y_i^{\ep}(z)Y_j^{\ep'}(w)&=&:Y_i^{\ep}(z)Y_j^{\ep'}(w):
(1-q^{-(\ep+\ep')/2}w/z)^{\ep\ep'(\a_i|\a_j)}_{q^2}\nonumber\\
&& \qquad \cdot z^{\ep\ep'(\a_i|\a_j)},
\label{OPE:1}\\
U_i(z)U_j(w)&=&:U_i(z)U_j(w):
(1-w/z)^{(\a_i|\a_j)}_{q^2}z^{(\a_i|\a_j)},  \label{OPE:2}\\
U_i(z)U_j^*(w)&=&:U_i(z)U_j^*(w):
(1-w/z)^{-(\a_i|\a_j)}_{q^2}z^{-(\a_i|\a_j)}. \label{OPE:3}
\end{eqnarray}
and the OPE's among $U^*_i(z)$ are the same as (\ref{OPE:2}).

\medskip
\noindent{\bf Proof of relation (\ref{E:R7})}.
There are four cases to be considered: $(\a_i|\a_j)=-1/2$, $(\a_i|\a_j)=-1$,
$(\a_i|\a_j)=1$, and $(\a_i|\a_j)=2$. Note that the construction of
$X^{\pm}_n(z)$ implies immediately that the last case holds because they
are the same as the simply laced cases.
Since the other verifications are similar,
we just show some computations in the following.

First let us consider $(\a_i|\a_j)=-1$:
\begin{eqnarray*}
 &&(z-qw)X^-_{n-1}(z)X^-_{n}(w)\\
 &=&:X^-_{n-1}(z)X^-_{n}(w):   \qquad \mbox{using (\ref{OPE:1})} \\
 &=&(qz-w)X^-_{n-1}(w)X^-_{n}(z)
\end{eqnarray*}

For simplicity in the following we will usually write
$$(z-w)^a_{q^2}=(1-w/z)^a_{q^2}z^a, $$
which is considered as a power series in $w/z$.

For $i=1, \cdots, n-2$ we have that
\begin{eqnarray*}
&&X^+_i(z)X^+_{i+1}(w)\\
&&=\vep_i\vep_{i+1}
:Y_i^+(z)Y_{i+1}^+(w):(z-q^{-1}w)_{q^2}^{-1/2} \\
 && \left(:U_i(zq^{1/2})U_{i+1}(wq^{1/2}):(z-w)^{-1/2}q^{1/4} \right.\\
&& +(-1)^{2a_{i+1}(0)}:U_i(zq^{1/2})U_{i+1}^*(wq^{-1/2}):
(z-qw)^{1/2}q^{-1/4}\\
 &&  +(-1)^{2a_{i}(0)+1}:U_i^*(zq^{-1/2})U_{i+1}(wq^{1/2}):
(z-q^{-1}w)^{1/2}q^{1/4}\\
&&  \left.+(-1)^{2a_i(0)+2a_{i+1}(0)}
:U_i^*(zq^{-1/2})U_{i+1}^*(wq^{-1/2}):(z-w)^{-1/2}
q^{-1/4}\right)
\end{eqnarray*}
Using $(z-w)^{-1/2}(z-q^{-1}w)^{-1/2}=(z-q^{-1/2}w)^{-1}$
we have
\begin{eqnarray*}
X^+_i(z)X^+_{i+1}(w)&=&
\vep_i\vep_{i+1}
:Y_i^+Y_{i+1}^+:\left(:U_iU_{i+1}(w-q^{-1/2}w)^{-1}q^{1/4}\right.\\
&&+(-1)^{2a_{i}(0)}:U_{i+1}U_{i}^*:\frac{z-q^{1/2}w}{z-q^{-1/2}}q^{-1/4}
\\
&&+(-1)^{1+2a_{i+1}(0)}:U_iU^*_{i+1}:q^{1/4}\\
&&\left.+(-1)^{2a_i(0)+2a_{i+1}(0)}
:U_i^*U_{i+1}^*(w-q^{-1/2}z)^{-1}q^{-1/4}\right)             ,
\end{eqnarray*}
>from which it follows that
\begin{eqnarray*}
&&(z-q^{-1/2}w)X_i^+(z)X_{i+1}^+(w)+(w-q^{-1/2}z)X_{i+1}^+(z)X_{i}^+(w)\\
&=&0.
\end{eqnarray*}

\medskip
\noindent{\bf Proof of relation (\ref{E:R8})}.
Again we need to consider the following four
cases: $(\a_i|\a_j)=-1/2$, $(\a_i|\a_j)=-1$,
$(\a_i|\a_j)=1$, and $(\a_i|\a_j)=2$. We only give the proof for
$i=j=1, \cdots, n-1$, since the other cases are either immediate or
similar to our previous considerations.

>From (\ref{OPE:1}-\ref{OPE:3}) it follows that
\begin{eqnarray*}
&&X^+_i(z)X_i^-(w)=\vep_i^2:Y_i^+(z)Y_i^+(w):\\
&&\left(:U_i(q^{1/2}z)U_i(q^{1/2}w):q^{1/2}\right. \\
&&+(-1)^{2a_i(0)}:U_i(q^{1/2}z)U_i^*(q^{-1/2}w):\frac 1{(q^{1/2}z-q^{-1/2}w)
(z-w)} \\
&&+(-1)^{2a_i(0)}:U_i^*(q^{-1/2}z)U_i(q^{1/2}w):\frac 1{(q^{-1/2}z-q^{1/2}w)
(z-w)} \\
&&+\left.(-1)^{4a_i(0)}:U_i^*(q^{-1/2}z)U_i^*(q^{-1/2}w):q^{-1/2}\right)
\end{eqnarray*}

Then we have
\begin{eqnarray*}
&&[X^+_{i}(z), X^-_{i}(w)]\\
&&=\vep_i^2:U_iU_i^*Y_i^+Y_i^-:\left(\frac1{(q^{1/2}z-q^{-1/2}w)
(z-w)}\right.\\
&&\qquad\quad \left. -\frac 1{(q^{-1/2}w-q^{1/2}z)
(w-z)}\right)\\
&&+\vep_i^2:U_i^*U_iY_i^+Y_i^-:\left(\frac1{(q^{-1/2}z-q^{1/2}w)
(z-w)}\right.\\
&&\qquad\quad \left.-\frac 1{(q^{1/2}w-q^{-1/2}z)
(w-z)}\right)\\
&&=\vep_i^2:U_iU_i^*Y_i^+Y_i^-:\frac1{(q^{1/2}-q^{-1/2})zw}\left(\delta(\frac
wz)-
\delta(q^{-1}\frac wz)\right)\\
&&+\vep_i^2:U_i^*U_iY_i^+Y_i^-:\frac1{(q^{1/2}-q^{-1/2})zw}\left(\delta(q\frac
wz)-
\delta(\frac wz)\right)\\
&&=\frac1{(q_i-q_i^{-1})zw}
                 \left(
                 \psi_i(wq^{1/2})
                 \delta(\frac{wq}{z})
                 -\varphi_i(wq^{-1/2})
                  \delta(\frac{wq^{-1}}z)
                 \right) .
\end{eqnarray*}

\medskip
\noindent{\bf Proof of Serre relations (\ref{E:R9})}.
For $i=1, \cdots, n-1$, let us write the operator $X^{\pm}_i(z)$ as a
sum of two terms:
\begin{eqnarray*}
X^{\pm}_i(z)&=&\sum_{\ep=\pm}Y_i^{\pm}(z)U_i^{\ep}(zq^{\pm\ep/2})
(-1)^{(1-\ep)a_i(0)}\\
&=& X^{\pm}_{i+}(z)+X^{\pm}_{i-}(z)
\end{eqnarray*}
where we identify $U^+_i(z)=U_i(z)$ and $U^-_i(z)=U^*_i(z)$.

>From the OPE (\ref{OPE:1}-\ref{OPE:3}) it follows that
\begin{eqnarray*}
&&X_{i\ep_1}^+(z_1)X_{i\ep_2}^+(z_2)X^+_{i+1,\ep}(w)=
:X_{i\ep_1}^+(z_1)X_{i\ep_2}^+(z_2)X^+_{i+1,\ep}(w):(-1)^{(\ep_1-\ep_2)/2}\\
&&\qquad \times
(z_1-q^{-1}z_2)(q^{\ep_1/2}z_1-q^{\ep_2/2}z_2)^{\ep_1\ep_2}
q^{-\ep/2}  \\
&&\qquad \times \frac{(z_1-q^{\ep/2}w)^{\frac{|\ep-\ep_1|}2}}{z_1-q^{-1/2}w}
\frac{(z_2-q^{\ep/2}w)^{\frac{|\ep-\ep_2|}2}}{z_2-q^{-1/2}w}
\end{eqnarray*}
where we include the sign factor $(-1)^{(1-\ep)\a_i(0)}$ and $\vep_i$ in
the normal ordered product. Similar normal product computation gives that

\begin{eqnarray}\label{E:S1}
&&X_{i\ep_1}^+(z_1)X_{i\ep_2}^+(z_2)X^+_{i+1,\ep}(w)-(q^{1/2}+q^{-1/2})
X_{i\ep_1}^+(z_1)X^+_{i+1,\ep}(w)X_{i\ep_2}^+(z_2)\nonumber\\
&&\quad +X^+_{i+1,\ep}(w)X_{i\ep_1}^+(z_1)X_{i\ep_2}^+(z_2)  \nonumber\\
&=&:X_{i\ep_1}^+(z_1)X_{i\ep_2}^+(z_2)X^+_{i+1,\ep}(w):
(z_1-q^{-1}z_2)(q^{\ep_1/2}z_1-q^{\ep_2/2}z_2)^{\ep_1\ep_2}q^{-\ep/2}
\nonumber\\
&& \quad \cdot \left( (-1)^{(\ep_1-\ep_2)/2}
\frac{(z_1-q^{\ep/2}w)^{\frac{|\ep-\ep_1|}2}
(z_2-q^{\ep/2}w)^{\frac{|\ep-\ep_2|}2}}{(z_1-q^{-1/2}w)(z_2-q^{-1/2}w)}\right.
\nonumber\\
&&\quad\quad +[2]_{q^{1/2}}(-1)^{(\ep_1-\ep_2)/2}
\frac{(z_1-q^{\ep/2}w)^{\frac{|\ep-\ep_1|}2}
(q^{\ep/2}w-z_2)^{\frac{|\ep-\ep_2|}2}}{(z_1-q^{-1/2}w)(w-q^{-1/2}z_2)}
     \nonumber\\
&&\qquad\quad \left. + \frac{(q^{\ep/2}w-z_1)^{\frac{|\ep-\ep_1|}2}
(q^{\ep/2}w-z_2)^{\frac{|\ep-\ep_2|}2}}{(w-q^{-1/2}z_1)(w-q^{-1/2}z_2)}
  \right)
\end{eqnarray}
where each term in the parentheses corresponds to the OPE's for the
three normal products, and we have also used the relation:
$$:X_{i\ep_1}^+(z_1)X_{i\ep_2}^+(z_2)X^+_{i+1,\ep}(w):
=-:X_{i\ep_1}^+(z_1)X^+_{i+1,\ep}(w)X_{i\ep_2}^+(z_2):.$$

We first claim that for $\ep_1\neq \ep_2$ we have
\begin{eqnarray}\label{E:S2}
&&X_{i\ep_1}^+(z_1)X_{i\ep_2}^+(z_2)X^+_{i+1,\ep}(w)-(q^{1/2}+q^{-1/2})
X_{i\ep_1}^+(z_1)X^+_{i+1,\ep}(w)X_{i\ep_2}^+(z_2)\nonumber\\
&&\quad +X^+_{i+1,\ep}(w)X_{i\ep_1}^+(z_1)X_{i\ep_2}^+(z_2)
+(z_1\leftrightarrow z_2, \ep_1\leftrightarrow \ep_2)=0
\end{eqnarray}

The claim is verified by checking 4 cases for $\ep, \ep_i$, which are all
similar and relied upon the following important identity \cite{kn:JKM}:
\begin{eqnarray}\label{identity1}
&&(z_1-aw)(z_2-aw)
+(a+a^{-1})(z_1-aw)(w-az_2)\nonumber\\
&&\qquad\qquad +(w-az_1)(w-az_2) \nonumber\\
&=&(a^{-1}-a)w(z_1-a^2z_2)
\end{eqnarray}
for any $a\in {\bf C}$.

In fact for $\ep=1, \ep_1=-\ep_2=1$, the parentheses in (\ref{E:S1})
is simplified to the following expression times
$q^{-1/2}\prod_i(z_i-q^{-1/2}w)^{-1}
\cdot (w-q^{-1/2}z_2)^{-1}$.

\begin{eqnarray*}
&&q^{-1/2}\left((z_1-q^{1/2}w)(z_2-q^{1/2}w)+[2]_i(z_1-q^{1/2}w)(w-q^{1/2}z_2)
\right.
\\
&&\qquad \left. +(w-q^{1/2}z_1)(w-q^{1/2}z_2)\right)\\
&&=(q^{-1/2}-q)w(z_1-q^{-1}z_2).
\end{eqnarray*}

Under the symmetry $(z_1, \ep_1)\leftrightarrow (z_2, \ep_2)$
it follows that the claim holds due to
$$w(z_1-q^{-1}z_2)+w(q^{-1}z_2-z_1)=0.$$

We now turn to the other 4 cases with $\ep_1=\ep_2$. The 4 cases are also
similar. Take the case $\ep=-\ep_1=-\ep_2=1$ for example. Using the identity
(\ref{identity1}) again to simplify the parentheses in (\ref{E:S1}), the
contraction function in the Serre relation becomes
\begin{eqnarray*}
&&q^{-1}(z_1-q^{-1}z_2)(z_1-z_2)
\left(\frac{(z_1-q^{1/2}w)(z_2-q^{1/2}w)}{(z_1-q^{-1/2}w)(z_2-q^{-1/2}w)}
\right.\\
&&\qquad\quad \left.-[2]_i\frac{z_1-q^{1/2}w}{z_1-q^{-1/2}w}q^{1/2}+q\right)\\
&=& \frac{q^{-1}(q^{-1/2}-q^{1/2})w
(z_1-z_2)(z_1-q^{-1}z_2)(z_1-qz_2)}{(z_1-q^{-1/2}w)(z_2-q^{-1/2}w)}
\end{eqnarray*}
which is anti-symmetric under $(z_1\leftrightarrow z_2)$, hence the sub-Serre
relation is proved in this case. That is,
\begin{eqnarray*}
&&X_{i\ep_1}^+(z_1)X_{i\ep_1}^+(z_2)X^+_{i+1,\ep}(w)-(q^{1/2}+q^{-1/2})
X_{i\ep_1}^+(z_1)X^+_{i+1,\ep}(w)X_{i\ep_1}^+(z_2)\\
&&\quad +X^+_{i+1,\ep}(w)X_{i\ep_1}^+(z_1)X_{i\ep_1}^+(z_2)
+(z_1\leftrightarrow z_2)=0
\end{eqnarray*}
Combining this sub-Serre relation with (\ref{E:S2}) we prove the Serre relation
for $A_{i, i+1}=A_{i+1. i}=-1$.

We remark that the case $A_{n-1, n}=-1$ is easily proved by using the
identity (\ref{identity1}) with $a=q$.

Finally let's show the fourth order Serre relation with $A_{n, n-1}=-2$.
\begin{eqnarray*}
&&Sym_{z_1, z_2, z_3}
(X^+_{n-1}(z_1)X^+_{n-1}(z_2)X^+_{n-1}(z_3)X_{n}^+(w)\\
&&\qquad -[3]_{q^{1/2}}
X^+_{n-1}(z_1)X^+_{n-1}(z_2)X_{n}^+(w)X^+_{n-1}(z_3)\\
&&\qquad  +[3]_{q^{1/2}}
X^+_{n-1}(z_1)X_{n}^+(w)X^+_{n-1}(z_2)X^+_{n-1}(z_3)\\
&&\qquad -X_{n}^+(w)X^+_{n-1}(z_1)X^+_{n-1}(z_2)X^+_{n-1}(z_3))=0
\end{eqnarray*}

First we have
\begin{eqnarray*}
&&X^+_{n-1, \ep_1}(z_1)X^+_{n-1, \ep_2}(z_2)X^+_{n-1, \ep_3}(z_3)
X_{n, \ep}^+(w)\\
&=&:X^+_{n-1, \ep_1}(z_1)\cdots:
\frac{\prod_{1\leq i<j\leq 3}(z_i-q^{-1}z_j)
(q^{\ep_i/2}z_i-q^{\ep_j/2}z_j)^{\ep_i\ep_j}}
{\prod_{i=1}^3
(z_i-q^{-1}w)(q^{\ep_i/2}z_i-q^{\ep/2}w)^{\ep_i\ep}}
\end{eqnarray*}

Pulling out the common normal product we have
\begin{eqnarray*}
&&X^+_{n-1, \ep_1}(z_1)X^+_{n-1, \ep_2}(z_2)X^+_{n-1, \ep_3}(z_3)
X_{n, \ep}^+(w)\\
&&\ \ -[3]_{q^{1/2}}X^+_{n-1, \ep_1}(z_1)X^+_{n-1, \ep_2}(z_2)
X_{n, \ep}^+(w)X^+_{n-1, \ep_3}(z_3)
\\
&&\ \ +[3]_{q^{1/2}}X^+_{n-1, \ep_1}(z_1)X_{n, \ep}^+(w)
X^+_{n-1, \ep_2}(z_2)X^+_{n-1, \ep_3}(z_3)\\
&&\ \ -X_{n, \ep}^+(w)X^+_{n-1, \ep_1}(z_1)X^+_{n-1, \ep_2}(z_2)
X^+_{n-1, \ep_3}(z_3)X_{n, \ep}^+(w)
\\
&&=:X^+_{n-1, \ep_1}(z_1)\cdots:
\frac{\prod_{1\leq i<j\leq 3}q^{-3}(z_i-q^{-1}z_j)
(q^{\ep_i/2}z_i-q^{\ep_j/2}z_j)^{\ep_i\ep_j}}
{\prod_{i=1}^3
(z_i-q^{-1}w)(w-q^{-1}z_i)(q^{\ep_i/2}z_i-q^{\ep/2}w)^{\ep_i\ep}}\\
&&\cdot \left((z_1-qw)(z_2-qw)(z_3-qw)+[3]_{q^{1/2}}
(z_1-qw)(z_2-qw)(w-qz_3)\right.\\
&&\left. +[3]_{q^{1/2}}(z_1-qw)(w-qz_2)(w-qz_3)+(w-qz_1)(w-qz_2)(w-qz_3)\right)
\end{eqnarray*}

The Serre relation is then equivalent to the following combinatorial identity
(cf. \cite{kn:J5}):
\begin{eqnarray}\label{identity2}
&&\sum_{\sigma\in S_3}(-1)^{l(\sigma)}\sigma .
\left\{(z_1-qw)(z_2-qw)(z_3-qw)\right. \nonumber\\
&&\qquad +[3]_{q^{1/2}}
(z_1-qw)(z_2-qw)(w-qz_3)\nonumber\\
&&\qquad +[3]_{q^{1/2}}(z_1-qw)(w-qz_2)(w-qz_3) \nonumber\\
&&\qquad \left. +(w-qz_1)(w-qz_2)(w-qz_3)\right\}
\prod_{1\leq i<j\leq 3}(z_i-q^{-1}z_j)=0
\end{eqnarray}
where the symmetric group $S_3$ acts on the ring of functions in $z_i$
($i=1, 2, 3$) in the natural way: $\sigma.z_i=z_{\sigma(i)}$.
Note that the expression in the parenthesis can be simplified to
$$
(q^{-1}-q)\left(w^2(z_1-(q+q^{-1})z_2+q^3z_3)+
w(z_1z_2-(q+q^{-1})z_1z_3+q^3z_2z_3)\right)
$$
Hence the identity (\ref{identity2}) is equivalent to the following.
\begin{equation}\label{identity3}
\sum_{\sigma\in S_3}sgn(\sigma) \sigma.(z_1-(q+q^2)z_2+q^3z_3)
\prod_{i<j}(qz_i-z_j)=0
\end{equation}
which was already proved in \cite{kn:JKM}. The proof of the Serre relation
is thus completed.

\noindent{\bf Highest weight vectors}.
To calculate the highest weight vectors we need the exact isomorphism
\cite{kn:J4}
between Drinfeld realizations and the Drinfeld-Jimbo definition of quantum
affine algebras. From \cite{kn:J4} we have
\begin{eqnarray*}
&&e_0=[X_1^-(0), \cdots, X_n^-(0), \\
&&\qquad X_{n-1}^-(0), \cdots,
X_1^-(1)]_{q^{-1/2}, \cdots, q^{-1/2}, q^{-1}, q^{-1/2}, \cdots, q^{-1/2}, 1}
\gamma K_{\theta}^{-1}\\
&&e_i=
X_i^+(0)
\end{eqnarray*}
where $K_{\theta}=K_1^2\cdots K_{n-1}^2K_n$. The $q$-multibracket is defined
inductively by
\begin{eqnarray*}
&&[a_1, a_2]_v=a_1a_2-va_2a_1\\
&&[a_1, a_2, \cdots, a_n]_{v_1, \cdots, v_{n-1}}=
[a_1, [a_2, \cdots, a_n]_{v_1, \cdots, v_{n-2}}]_{v_{n-1}}
\end{eqnarray*}

Standard calculation of contour integrals of vertex operators will immediately
give the following result.

\begin{lemma} (i) Let $\lambda\in P$ and $\tla\in \tilde{P}$ be dominant
with $(\la, \a_i)=(\tla, \ta_i)$, then for any $j=1, \cdots, n$ and $n\geq 0$
we have
$$
X_j^+(n).e^{\la}e^{\tla}=0.
$$
(ii) For the fundamental weights $\la_i\in P$, $\tla_i\in \tilde{P}$
we have for $j=1, \cdots, n$
$$
X_j^-(0)e^{\la}e^{\tla}=-\ep(\a_i, \lambda_i)e^{\la_i-\a_i}e^{\tla_i-\ta_i}
\delta_{ij}, \ \ X_j^-(1)e^{\la_i}e^{\tla_i}=0.
$$
\end{lemma}

Using this lemma and the isomorphism of $e_0$ and $e_i$ it is easy to see that
$e^{\la_i}e^{\tla_i}$ ($i=0, \cdots, n$) are indeed
 highest weight vectors contained in the module $\cal V$.

\medskip

\centerline{\bf Acknowledgments}
The first author (NJ) is  supported in part by NSA
grants MDA 904-96-1-0087 and MDA 904-97-1-0062. The second author
(YK) is supported by the JSPS Research Fellowships for Young Scientists.
The third author (KCM)
is supported in part by NSA grant MDA 904-96-1-0013.
The second author (YK) thanks the mathematics department at N. C. State
University for the hospitality during his visit when part of this work was
done.

\end{document}